\documentclass[12pt]{article}
\usepackage[utf8]{inputenc}
\usepackage{mathrsfs}
\usepackage{fullpage}
\usepackage{amssymb}
\usepackage{amsmath}
\usepackage{epsfig}
\usepackage{algorithm}
\usepackage{algorithmic}
\usepackage{latexsym}
\usepackage{amsmath}
\usepackage{amsfonts}
\usepackage{graphicx}
\usepackage{authblk}
\usepackage[normalem]{ulem}
\usepackage{array}
\usepackage{colortbl}
\usepackage{color}\usepackage[dvipsnames]{xcolor}
\usepackage{tikz}\usetikzlibrary{patterns}
\usepackage{array,multirow}
\usepackage[
  bookmarks=false, 
  colorlinks,
  citecolor=brown!70!black,
  linkcolor=brown!80!black,
  urlcolor=blue!70!black,
]{hyperref}
\usepackage{relsize,exscale}

\newtheorem{thm}{Theorem}

\newcommand{\comm}[1]{}
 	\definecolor{lightlightgray}{rgb}{0.93, 0.93, 0.93}
 		\definecolor{llightgray}{rgb}{0.87, 0.87, 0.87}

\newcolumntype{C}[1]{>{\centering\arraybackslash }b{#1}}

\def\A{0.5cm}

 \def\C{2.5cm}
 
 \def\E{4.5cm}

 \def\G{6.5cm}
 
 \def\I{8.5cm}

 \def\K{10.5cm}
 
\def\M{12.5cm}
 
 \def\O{14.5cm}\def\Q{16.5cm}\def\S{18.5cm}
\def\U{20.5cm}\def\W{22.5cm}\def\Y{24.5cm}\def\ZZ{26.5cm}\def\Za{28.5cm}\def\Zb{30.5cm}\def\Zc{32.5cm}\def\Zd{34.5cm}\def\Ze{36.5cm}

\def\sizePoint{2.6pt}

\def\sizeCircle{5pt}

\newcommand{\point}[2]{\fill (canvas cs:x=#1,y=#2) circle (\sizePoint); }

\newcommand{\pointCercle}[2]{\point{#1}{#2} \draw (canvas cs:x=#1,y=#2) circle (\sizeCircle); }
\newcommand\oeis[1]{\href{https://oeis.org/#1}{#1}}

\title{Enumeration of partial \L{}ukasiewicz paths}

\author[1]{Jean-Luc Baril}
\author[2]{Helmut Prodinger}
\affil[1]{\rm LIB, Universit\'e de Bourgogne Franche-Comt\'e \protect\\
  B.P. 47 870, 21078 Dijon Cedex France\protect\\
   {\tt E-mail: barjl@u-bourgogne.fr
   }
}
\affil[2]{\rm Department of Mathematical Sciences, Stellenbosch University\protect\\
7602 Stellenbosch, South Africa\protect\\
   {\tt E-mail: hproding@sun.ac.za
   }
}

\date{April 2022}

\begin{document}

\maketitle
\begin{abstract} \L{}ukasiewicz  paths are lattice paths in $\Bbb{N}^2$ starting at the origin, ending on the $x$-axis, and consisting of steps in the set $\{(1,k), k\geq -1\}$. We give generating function and exact value for the number of $n$-length prefixes (resp.\ suffixes) of these paths  ending at height $k\geq 0$ with a given type of step. We make a similar study for prefixes of height at most $t\geq 0$.
	Using the explicit forms for the paths of bounded height, we evaluate the average height asymptotically. For fixed $k$ and $n\to\infty$, this quantity behaves as
	$\sqrt{\pi n}$.
	Finally we study (in the same way) prefixes of alternate \L{}ukasiewicz  paths, i.e., \L{}ukasiewicz  paths that do contain two consecutive steps with the same direction.
\end{abstract}
\section{Introduction}

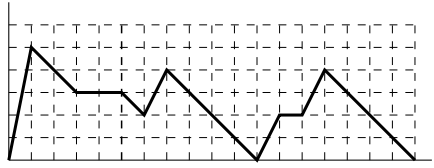
\begin{figure}[h]
 \begin{center}
        \begin{tikzpicture}[scale=0.15]
            \draw (\A,\A)-- (38,\A);
             \draw[dashed,line width=0.1mm] (\A,\E)-- (\Ze,\E);
              \draw[dashed,line width=0.1mm] (\A,\C)-- (\Ze,\C);
               \draw[dashed,line width=0.1mm] (\A,\G)-- (\Ze,\G);
               \draw[dashed,line width=0.1mm] (\A,\I)-- (\Ze,\I);
              \draw[dashed,line width=0.1mm] (\A,\K)-- (\Ze,\K);
               \draw[dashed,line width=0.1mm] (\A,\M)-- (\Ze,\M);
            \draw (\A,\A) -- (\A,\O);
             \draw[dashed,line width=0.1mm] (\C,\A) -- (\C,\M);\draw[dashed,line width=0.1mm] (\E,\A) -- (\E,\M);\draw[dashed,line width=0.1mm] (\G,\A) -- (\G,\M);
             \draw[dashed,line width=0.1mm] (\I,\A) -- (\I,\M);\draw[dashed,line width=0.2mm] (\K,\A) -- (\K,\M);\draw[dashed,line width=0.1mm] (\M,\A) -- (\M,\M);
             \draw[dashed,line width=0.1mm] (\O,\A) -- (\O,\M);\draw[dashed,line width=0.1mm] (\Q,\A) -- (\Q,\M);\draw[dashed,line width=0.1mm] (\S,\A) -- (\S,\M);
             \draw[dashed,line width=0.1mm] (\U,\A) -- (\U,\M);\draw[dashed,line width=0.1mm] (\W,\A) -- (\W,\M);\draw[dashed,line width=0.1mm] (\Y,\A) -- (\Y,\M);
             \draw[dashed,line width=0.1mm] (\ZZ,\A) -- (\ZZ,\M);
             \draw[dashed,line width=0.1mm] (\Za,\A) -- (\Za,\M);
             \draw[dashed,line width=0.1mm] (\Zb,\A) -- (\Zb,\M);
             \draw[dashed,line width=0.1mm] (\Zc,\A) -- (\Zc,\M);
             \draw[dashed,line width=0.1mm] (\Zd,\A) -- (\Zd,\M);
             \draw[dashed,line width=0.1mm] (\Ze,\A) -- (\Ze,\M);
            \draw[solid,line width=0.4mm] (\A,\A)--(\C,\K) -- (\E,\I) -- (\G,\G) -- (\I,\G) --(\K,\G)-- (\M,\E) -- (\O,\I) -- (\Q,\G) -- (\S,\E) -- (\U,\C)  -- (\W,\A) --(\Y,\E) -- (\ZZ,\E) --(\Za,\I) -- (\Ze, \A);

         \end{tikzpicture}

               \end{center}
         \caption{ A \L{}ukasiewicz path of length $18$: $U_5DDFFDU_2DDDDU_2FU_2DDDD$.}
         \label{fig1}
\end{figure}

A {\it  \L{}ukasiewicz path} of length $n\geq 0$ is a lattice path in $\Bbb N^2$ starting at the origin $(0,0)$, ending on the $x$-axis, consisting of $n$ steps lying in $S=\{(1,k), k\geq -1\}$.  We denote by $\epsilon$ the empty path, {\it i.e.}, the path of length zero. These paths constitute a natural generalization of Dyck and  Motzkin paths (see \cite{Deu,Don}), which are made using  steps into the sets $\{(1,1),(1,-1)\}$ and $\{(1,1),(1,0), (1,-1)\}$, respectively. We refer to \cite{Bar0,Ges,Ran,Sta,Var,Vie} for some combinatorial studies on \L{}ukasiewicz paths. Let $\mathcal{L}_n$, $n\geq 0$, be the set of \L{}ukasiewicz paths of length  $n$, and $\mathcal{L}=\bigcup_{n\geq 0} \mathcal{L}_n$. For convenience, we set $D=(1,-1)$, $F=(1,0)$, $U_k=(1,k)$ for $k\geq 1$. See Figure \ref{fig1} for an illustration of a  \L{}ukasiewicz path of length $18$. Note that \L{}ukasiewicz paths can be interpreted as an algebraic language  of words $w\in \{x_0,x_1,x_2, \ldots\}^\star$
   such that $\delta(w)=-1$ and $\delta(w')\geq 0$ for any proper prefix $w'$ of $w$ where $\delta$ is the map from $\{x_0,x_1,x_2, \ldots\}^\star$ to $\mathbb{Z}$ defined by $\delta(w_1w_2\ldots w_n)=\sum_{i=1}^n\delta(w_i)$
   with $\delta(x_i)=i-1$ (see \cite{Lot,Sch}).

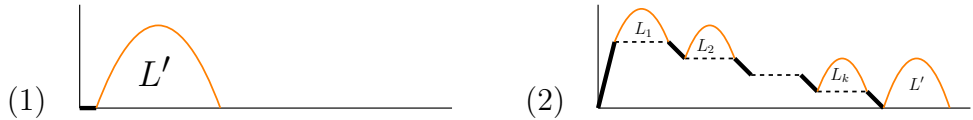
\begin{figure}[h]
\begin{center}(1)\quad\scalebox{0.55}{\begin{tikzpicture}[ultra thick]
 \draw[black, thick] (0,0)--(9,0); \draw[black, thick] (0,0)--(0,2.5);
  \draw[black, line width=3pt] (0,0)--(0.4,0);
  \draw[orange,very thick] (0.4,0) parabola bend (1.9,2) (3.4,0);
 \draw  (1.8,0.9) node {\huge $L'$};
 \end{tikzpicture}}
 \qquad (2)\quad\scalebox{0.55}{\begin{tikzpicture}[ultra thick]
 \draw[black, thick] (0,0)--(9,0); \draw[black, thick] (0,0)--(0,2.5);
  \draw[black, line width=3pt] (0,0)--(0.4,1.6);
  \draw[black, dashed, very thick] (0.4,1.6)--(1.7,1.6);
   \draw[black, line width=3pt] (1.7,1.6)--(2.1,1.2);
   \draw[black, dashed, very thick] (2.1,1.2)--(3.3,1.2);
   \draw[black, line width=3pt] (3.3,1.2)--(3.7,0.8);
   \draw[black, dashed, very thick] (3.7,0.8)--(4.9,0.8);
    \draw[black, line width=3pt] (4.9,0.8)--(5.3,0.4);
   \draw[black, dashed, very thick] (5.3,0.4)--(6.5,0.4);
   \draw[black, line width=3pt] (6.5,0.4)--(6.9,0);
 \draw[orange,very thick] (0.4,1.6) parabola bend (1.05,2.4) (1.7,1.6);
  \draw[orange,very thick] (2.1,1.2) parabola bend (2.7,2) (3.3,1.2);
  \draw[orange,very thick] (5.3,0.4) parabola bend (5.9,1.2) (6.5,0.4);
  \draw[orange,very thick] (6.9,0) parabola bend (7.7,1.2) (8.5,0);
 \draw  (1.1,1.9) node {$L_1$};\draw  (2.6,1.5) node { $L_2$};
 \draw  (5.85,0.8) node { $L_k$};\draw  (7.7,0.6) node { $L'$};
 \end{tikzpicture}}

\end{center}
\caption{The two forms  of  a non-empty \L{}ukasiewicz path.}
\label{fig2}\end{figure}

Any non-empty \L{}ukasiewicz path $L\in \mathcal{L}$ can be decomposed (see \cite{Fla}) into one of the two following forms: (1) $L=FL'$ with $L'\in\mathcal{L}$, or (2) $L=U_kL_1DL_2D\ldots L_kD L'$ with $k\geq 1$ and $L_1,L_2, \ldots, L_k, L'\in \mathcal{L}$ (see Figure \ref{fig2}).
Due to this decomposition, the generating function $L(z)=\sum_{n\geq 0}a_nz^n$ where $a_n$ is the cardinality of $\mathcal{L}_n$, satisfies the functional equation $L(z)=1+zL(z)+\sum_{k\geq 1}z^{k+1}L(z)^{k+1}$,  or equivalently,  $L(z)=\frac{1}{1-zL(z)}$. Then, $L(z)=\frac{1-\sqrt{1-4z}}{2z}$ and $a_n$ is the $n$-th Catalan number $a_n=\frac{1}{n+1}{2n \choose n}$ (see sequence \href{https://oeis.org/A000108}{A000108} in \cite{Sloa}).

In this paper, we provide enumerating results for several classes of partial \L{}ukasiewicz paths (prefixes and suffixes of \L{}ukasiewicz paths, partial alternate  \L{}ukasiewicz paths). More precisely, in Sections 2 and 3, we give generating functions and exact values for the number of $n$-length prefixes (resp.\ suffixes) of these paths ending at height $k\geq 0$ with a given type of step (down, up or horizontal step). In Section 4, We make a similar study for the paths of height at most $t\geq 0$.
Section 5 deals with the average height of various families of the \L{}ukasiewicz type.
In Section 6, we focus on partial alternate \L{}ukasiewicz  paths, i.e., \L{}ukasiewicz  paths that do contain two consecutive steps with the same direction.

All our \emph{explicit} formul\ae \ follow from the standard identity
$$[z^n]\biggl(\frac{1-\sqrt{1-4z}}{2z}\biggr)^k=\binom{2n-1+k}{n}-\binom{2n-1+k}{n}.$$


\section{Enumeration of Partial Lukasiewicz paths}

\begin{figure}[h]
\begin{center}\quad\scalebox{0.85}{\begin{tikzpicture}[ultra thick]
 \node at (0,2) {} edge [in=130,out=190,loop, thick,blue] ();
 \node at (2,2) {} edge [in=130,out=190,loop, thick,blue] ();
 \node at (4,2) {} edge [in=130,out=190,loop, thick,blue] ();
 \node at (6,2) {} edge [in=130,out=190,loop, thick,blue] ();
 \node at (8,2) {} edge [in=130,out=190,loop, thick,blue] ();
 \pointCercle{0cm}{4cm}
 \point{2cm}{4cm}\point{4cm}{4cm}\point{6cm}{4cm}\point{8cm}{4cm}
\point{0cm}{2cm}
\point{2cm}{2cm}\point{4cm}{2cm}\point{6cm}{2cm}\point{8cm}{2cm}
\point{0cm}{0cm}\point{2cm}{0cm}\point{4cm}{0cm}\point{6cm}{0cm}\point{8cm}{0cm}
  \draw[->,>=latex,black,thick] (0,4)--(1.92,4);
  \draw[->,>=latex,black,thick] (2,4)--(3.92,4);
  \draw[->,>=latex,black,thick] (4,4)--(5.92,4);
  \draw[->,>=latex,black,thick] (6,4)--(7.92,4);
    \draw[black,dashed,thick] (8.5,4)--(9,4);
    \draw[black,dashed,thick] (8.5,2)--(9,2);
    \draw[black,dashed,thick] (8.5,0)--(9,0);
  \draw[->,>=latex,thick] (0,4) to[bend left=20] (3.99,4);
  \draw[->,>=latex,thick] (0,4) to[bend left=20] (5.99,4);
  \draw[->,>=latex,thick] (0,4) to[bend left=20] (7.99,4);
  \draw[->,>=latex,thick] (2,4) to[bend left=20] (3.99,4);
  \draw[->,>=latex,thick] (2,4) to[bend left=20] (5.99,4);
  \draw[->,>=latex,thick] (2,4) to[bend left=20] (7.99,4);
  \draw[->,>=latex,thick] (4,4) to[bend left=20] (7.99,4);
  \draw[->,>=latex,thick] (0,2) to[bend left=20] (1.92,3.92);
    \draw[->,>=latex,thick] (0,2) to[bend left=15] (3.92,3.92);
     \draw[->,>=latex,thick] (0,2) to[bend left=10] (5.92,3.92);
     \draw[->,>=latex,thick] (0,2) to[bend left=5] (7.92,3.92);
     \draw[->,>=latex,thick] (0,0) to[bend left=20] (1.92,3.92);
    \draw[->,>=latex,thick] (0,0) to[bend left=15] (3.92,3.92);
     \draw[->,>=latex,thick] (0,0) to[bend left=10] (5.92,3.92);
     \draw[->,>=latex,thick] (0,0) to[bend left=9] (7.92,3.92);
     \draw[->,>=latex,red,thick] (2,2) to[bend right=5] (0.05,0.05);
      \draw[->,>=latex,red,thick] (4,2) to[bend right=20] (2.01,0.05);
     \draw[->,>=latex,red,thick] (6,2) to[bend right=20] (4.01,0.05);
     \draw[->,>=latex,red,thick] (8,2) to[bend right=20] (6.01,0.05);
     \draw[->,>=latex,red,thick] (2,4)-- (0.02,0.05);
     \draw[->,>=latex,red,thick] (4,4)-- (2.02,0.05);
     \draw[->,>=latex,red,thick] (6,4)-- (4.02,0.05);
     \draw[->,>=latex,red,thick] (8,4)-- (6.02,0.05);
      \draw[thick] (2,0)-- (2.4,0.57);
      \draw[->,>=latex,thick] (2,0)-- (2.6,0.6);
       \draw[thick] (2,0)-- (2.6,0.4);
       \draw[thick] (2,2)-- (2.4,2.57);
      \draw[->,>=latex,thick] (2,2)-- (2.6,2.6);
       \draw[thick] (2,2)-- (2.6,2.4);
        \draw[thick] (4,0)-- (4.4,0.57);
      \draw[->,>=latex,thick] (4,0)-- (4.6,0.6);
       \draw[thick] (4,0)-- (4.6,0.4);
       \draw[thick] (4,2)-- (4.4,2.57);
      \draw[->,>=latex,thick] (4,2)-- (4.6,2.6);
       \draw[thick] (4,2)-- (4.6,2.4);
        \draw[thick] (6,0)-- (6.4,0.57);
      \draw[->,>=latex,thick] (6,0)-- (6.6,0.6);
       \draw[thick] (6,0)-- (6.6,0.4);
       \draw[thick] (6,2)-- (6.4,2.57);
      \draw[->,>=latex,thick] (6,2)-- (6.6,2.6);
       \draw[thick] (6,2)-- (6.6,2.4);
        \draw[thick] (8,0)-- (8.4,0.57);
      \draw[->,>=latex,thick] (8,0)-- (8.6,0.6);
       \draw[thick] (8,0)-- (8.6,0.4);
       \draw[thick] (4,4)-- (4.4,4.57);
      \draw[->,>=latex,thick] (4,4)-- (4.6,4.6);
       \draw[thick] (4,4)-- (4.6,4.4);
        \draw[thick] (6,4)-- (6.4,4.57);
      \draw[->,>=latex,thick] (6,4)-- (6.6,4.6);
       \draw[thick] (6,4)-- (6.6,4.4);
       \draw[thick] (8,2)-- (8.4,2.57);
      \draw[->,>=latex,thick] (8,2)-- (8.6,2.6);
       \draw[thick] (8,2)-- (8.6,2.4);
       \draw[thick] (8,4)-- (8.4,4.57);
      \draw[->,>=latex,thick] (8,4)-- (8.6,4.6);
       \draw[thick] (8,4)-- (8.6,4.4);
\draw[->,>=latex,blue, thick] (0,0)--(0,1.92);
\draw[->,>=latex,blue, thick] (2,0)--(2,1.92);
\draw[->,>=latex,blue, thick] (4,0)--(4,1.92);
\draw[->,>=latex,blue, thick] (6,0)--(6,1.92);
\draw[->,>=latex,blue, thick] (8,0)--(8,1.92);
\draw[->,>=latex,blue, thick] (0,4)--(0,2.08);
\draw[->,>=latex,blue, thick] (2,4)--(2,2.08);
\draw[->,>=latex,blue, thick] (4,4)--(4,2.08);
\draw[->,>=latex,blue, thick] (6,4)--(6,2.08);
\draw[->,>=latex,blue, thick] (8,4)--(8,2.08);
  \draw[<-,>=latex,red,thick] (0.08,0)--(2,0);
  \draw[<-,>=latex,red,thick] (2.08,0)--(4,0);
  \draw[<-,>=latex,red,thick] (4.08,0)--(6,0);
  \draw[<-,>=latex,red,thick] (6.08,0)--(8,0);
 %
%
 \end{tikzpicture}}
\end{center}
\caption{The state diagram for the generation of partial \L{}ukasiewicz paths. Black (resp.\ red, resp.\ blue) arrows corresponds to up-steps (resp.\ down-steps, resp.\ horizontal steps).}
\label{fig3}\end{figure}
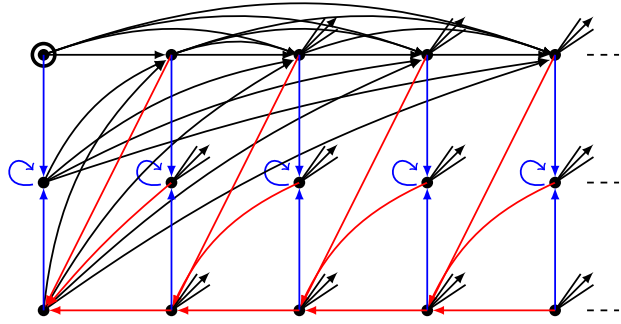

Partial  \L{}ukasiewicz paths of length $n$ (i.e. $n$-length prefixes of \L{}ukasiewicz paths) ending at height $k$ can be constructed trough the following state diagram (Figure~\ref{fig3}). We start at the first state of the first layer and end at the $(k+1)$-th state of one of the three layers, knowing that:  a move  ending at the $j$-th step of  the first layer (black arrows) and starting at the $i$-th state ($i<j$) of any layer corresponds to an up-step $U_{j-i}$, a move  ending at the $i$-th step of  the second layer (blue arrows) and starting at the $i$-th state of any layer corresponds to an horizontal step $F$, and a move  ending at the $i$-th step of  the third layer (red arrows) and starting at the $(i+1)$-th state of any layer corresponds to a down-step $D$.

For $k\geq 0$, we consider the generating function $f_k=f_k(z)$ (resp.\ $g_k=g_k(z)$, resp.\ $h_k=h_k(z)$) where the coefficient of $z^n$ in the series expansion is the number of  partial \L{}ukasiewicz paths of length $n$  ending at height $k$ with an up-step $U_k$, $k\geq 1$, (resp.\ with a down-step $D$, resp.\  with an horizontal-step $F$). Considering the state diagram in Figure~\ref{fig3}, $f_k$ (resp.\ $g_k$, resp.\ $h_k$) is the generating function in the variable $z$ marking the length of the paths ending on the $(k+1)$-th state of the first (resp.\ second, resp.\  third) layer.  So, we easily obtain the following equations:

\begin{equation}\begin{array}{l}
f_0=1,\mbox{ and } f_k=z\sum\limits_{\ell=0}^{k-1} f_\ell+z\sum\limits_{\ell=0}^{k-1} g_\ell +z\sum\limits_{\ell=0}^{k-1} h_\ell, \quad k\geq 1,\\
g_k=z f_{k+1}+zg_{k+1}+zh_{k+1},\quad k\geq 0,\\
h_k=zf_k+zg_k+zh_k,\quad k\geq 0.\\
\end{array}\label{equ}
\end{equation}

Now we introduce bivariate generating functions:
$$F(u,z)=\sum\limits_{k\geq 0} u^kf_k(z), \quad G(u,z)=\sum\limits_{k\geq 0} u^kg_k(z), \mbox{ and } H(u,z)=\sum\limits_{k\geq 0} u^kh_k(z).$$
For short, we also use the notations $F(u), G(u)$ and $H(u)$ for these functions.
Summing the recursions in (\ref{equ}), we have:
\begin{align*}
F(u)&=1+z\sum\limits_{k\geq 1}u^k\biggl(\sum\limits_{\ell=0}^{k-1} f_\ell+\sum\limits_{\ell=0}^{k-1} g_\ell +\sum\limits_{\ell=0}^{k-1} h_\ell\biggr)\\
&=1+z\sum\limits_{k\geq 0}\frac{u^{k+1}}{1-u}f_k+z\sum\limits_{k\geq 0}\frac{u^{k+1}}{1-u}g_k+z\sum\limits_{k\geq 0}\frac{u^{k+1}}{1-u}h_k\\
&=1+\frac{uz}{1-u}(F(u)+G(u)+H(u)),\\
G(u)&=z\sum\limits_{k\geq 0}u^k \Bigl(f_{k+1}+ g_{k+1} + h_{k+1}\Bigr)\\
&=\frac{z}{u}(F(u)+G(u)+H(u)-F(0)-G(0)-H(0)),\\
H(u)&=\frac{z}{1-z}(F(u)+G(u))
\end{align*}
where  $F(0)+G(0)+H(0)$ is the number of \L{}ukasiewicz paths of length $n$, i.e.,
  $$F(0)+G(0)+H(0) =L(z)=\frac{1-\sqrt{1-4z}}{2z}.$$

Solving these functional equations, we deduce
$$F(u)=1-z-{\frac {z \left( 1+\sqrt {1-4z} \right) }{2 u-1-\sqrt {1-4 z}
}}, \quad G(u)={\frac {\sqrt {1-4 z}+2 z-1}{2 u-1-\sqrt {1-4 z}}}, \mbox{ and }$$
$$H(u)=z+{\frac {z \left( \sqrt {1-4 z}-1 \right) }{2 u-1-\sqrt {1-4 z}}},$$
which implies that
\begin{align}f_k&=[ u^k] F(u)=\frac{2^kz}{(1+\sqrt{1-4z})^k}=z\biggl(\frac{1-\sqrt{1-4z}}{2z}\biggr)^k,\\
g_k&=[ u^k] G(u)=\frac{2^k(1-2z-\sqrt{1-4z})}{(1+\sqrt{1-4z})^{k+1}}\nonumber\\
&=\frac{(1-2z-\sqrt{1-4z})}{2}\biggl(\frac{1-\sqrt{1-4z}}{2z}\biggr)^{k+1}\nonumber\\
&=z\biggl(\frac{1-\sqrt{1-4z}}{2z}\biggr)^{k+2}-z\biggl(\frac{1-\sqrt{1-4z}}{2z}\biggr)^{k+1},\\
 \intertext{ and }
h_k&=[ u^k] H(u)=\frac{2^kz(1-\sqrt{1-4z})}{(1+\sqrt{1-4z})^{k+1}}
=\frac{z(1-\sqrt{1-4z})}{2}\biggl(\frac{1-\sqrt{1-4z}}{2z}\biggr)^{k+1}\nonumber\\
&=z^2\biggl(\frac{1-\sqrt{1-4z}}{2z}\biggr)^{k+2}.
\end{align}

\begin{thm} The bivariate generating function for the total number of partial \L{}ukasiewicz paths of length $n$ with respect to the height of the end point is given by
$$\mathit{Total}(z,u)=1+\frac{-1+\sqrt{1-4z}}{2u-1-\sqrt{1-4z}},$$ and we have
$$[ u^k] \mathit{Total}(z,u)=[k=0]+zL(z)^{k+2}.$$
Finally we have for $n\geq 1$, $$[ z^n][ u^k]\mathit{Total}(z,u)=\frac{k+2}{n+k+1}{2n+k-1\choose n-1},$$
$$[ z^n][ u^k]F(u)=\frac{k}{n+k-1}{2n+k-3\choose n-1},$$
$$[ z^n][ u^k]G(u)=\frac{k+3}{n+k+1}{2n+k-2\choose n-2},$$
$$[ z^n][ u^k]H(u)=\frac{k+2}{n+k}{2n+k-3\choose n-2}.$$
\end{thm}
Here are examples of  the series expansions of $[ u^k]\mathit{Total}(z,u)$ for $k=0,1,2,3$ (leading terms):

\noindent $\bullet \quad 1+z+2z^2+5z^3+14z^4+42z^5+132z^6+429z^7+1430z^8+4862z^9,$

\noindent $\bullet \quad z+3z^2+9z^3+28z^4+90z^5+297z^6+1001z^7+3432z^8+11934z^9,$

\noindent $\bullet \quad z+4z^2+14z^3+48z^4+165z^5+572z^6+2002z^7+7072z^8+25194z^9,$

\noindent $\bullet \quad z+5z^2+20z^3+75z^4+275z^5+1001z^6+3640z^7+13260z^8+48450z^9,$

\noindent which correspond respectively to \oeis{A000108}, \oeis{A000245}, \oeis{A002057}, and \oeis{A000344} in \cite{Sloa}.

According to Theorem~3.1 and Theorem~3.3 in \cite{Gud},  $[ z^n][ u^k]\mathit{Total}(z,u)$ counts also standard Young tableaux of shape $(n+2, n-k+1)$ (see \cite{Knu,You} for the definition of a standard Young tableau), and Dyck paths of semilength $n+k$ starting with at least $k$ up-steps and touching the $x$-axis somewhere between the two endpoints.


\section{Partial \L{}ukasiewicz paths from right to left}

In this section, we count partial \L{}ukasiewicz paths read from right to left, i.e., paths in $\Bbb{N}^2$ starting at the origin, consisting of steps $(1,k)$ for $k\leq 1$, and ending at a given height with a given type of step.  Of course, this study is completely equivalent to count suffixes of \L{}ukasiewicz paths starting at a given height with a given type of step. The notations $f_k$, $g_k$ and $h_k$ are  generating functions for the number of these paths (with respect to the length) ending at height $k$ with an up-step, down-step, horizontal step respectively.

Then we have,
\begin{align}
		f_0&=1\quad\mbox{and}\quad f_k=zf_{k-1}+zg_{k-1}+zh_{k-1},\quad k\geq 1,\nonumber\\
		g_k&=z\sum\limits_{\ell\geq k+1} f_\ell+z\sum\limits_{\ell\geq k+1} g_\ell +z\sum\limits_{\ell\geq k+1} h_\ell,\quad k\geq 0,\\
		h_k&=zf_k+zg_k+zh_k,\quad k\geq 0.\nonumber
\label{equ2}
\end{align}

Considering the bivariate generating functions
$$F(u)=\sum\limits_{k\geq 0} u^kf_k(z), \quad G(u)=\sum\limits_{k\geq 0} u^kg_k(z), \mbox{ and } H(u)=\sum\limits_{k\geq 0} u^kh_k(z),$$ and summing the recursions in (\ref{equ2}), we obtain:

\begin{align*}
F(u)&=1+z\sum\limits_{k\geq 1}u^k \left(f_{k-1}+ g_{k-1} + h_{k-1}\right)\\
&=1+zuF(u)+zuG(u)+zuH(u);\\
G(u)&=z\sum\limits_{k\geq 0}u^k\biggl(\sum\limits_{\ell\geq k+1} f_\ell+\sum\limits_{\ell\geq k+1} g_\ell +\sum\limits_{\ell\geq k+1} h_\ell\biggr)\\
&=z\sum\limits_{k\geq 1}\frac{1-u^{k}}{1-u}f_k+z\sum\limits_{k\geq 1}\frac{1-u^{k}}{1-u}g_k+z\sum\limits_{k\geq 1}\frac{1-u^{k}}{1-u}h_k\\
&=\frac{z}{1-u}(F(1)+G(1)+H(1)-F(u)-G(u)-H(u));\\
H(u)&=\frac{z}{1-z}(F(u)+G(u))\\
\end{align*}
with $$F(0)+G(0)+H(0)=L(z)=\frac{1-\sqrt{1-4z}}{2z}.$$
Moreover, we have $$F(1)+G(1)+H(1)=\frac{L(z)-1}{z}$$ because there is a bijection between all partial \L{}ukasiewicz paths read from right to left of length $n$ and \L{}ukasiewicz paths read from left to right of length $n+1$ (it suffices to remove the first step of every \L{}ukasiewicz path, and to read it from right to left).

Solving these functional equations, we deduce
$$F(u)=-\frac{1+\sqrt{1-4 z}}{2 z u -\sqrt{1-4 z}-1},\qquad G(u)= \frac{-1+\sqrt{1-4 z}+2 z}{2 z u -\sqrt{1-4 z}-1},$$
$$H(u)=-\frac{2 z}{2 z u -\sqrt{1-4 z}-1}$$
which implies that
\begin{align}f_k&= [ u^k] F(u)=\frac{2^kz^k}{(1+\sqrt{1-4z})^k},\\
g_k&=[ u^k] G(u)=\frac{2^kz^k(1-2z-\sqrt{1-4z})}{(1+\sqrt{1-4z})^{k+1}}, \mbox{ and }\\
h_k&=[ u^k] H(u)=\frac{2^{k+1}z^{k+1}}{(1+\sqrt{1-4z})^{k+1}}.
\end{align}

\begin{thm} The bivariate generating function for the total number of partial \L{}ukasiewicz paths of length $n$ (read from right to left) with respect to the height of the end point is given by
$$\mathit{Total}(z,u)=1+\frac{2}{1-2zu+\sqrt{1-4z}},$$ and we have
$$[ u^k] \mathit{Total}(z,u)=z^kL(z)^{k+1}.$$
Finally we have for $n\geq 1$, $$[ z^n][ u^k]Total(z,u)=\frac{k+1}{n+1}{2n-k\choose n},$$
\begin{gather*}
[ z^n][ u^k]F(u)=\frac{k}{n}{2n-k-1\choose n-1},\\
[ z^n][ u^k]G(u)=\frac{k+3}{n+1}{2n-k-2\choose n},\\
[ z^n] [ u^k ]H(u)=\frac{k+1}{n}{2n-k-2\choose n-1}.
\end{gather*}

\end{thm}
Here are examples of  the series expansions of $ [ u^k ]\mathit{Total}(z,u)$ for $k=0,1,2,3$ (leading terms):

\noindent $\bullet \quad 1+z+2z^2+5z^3+14z^4+42z^5+132z^6+429z^7+1430z^8+4862z^9,$

\noindent $\bullet \quad z+2z^2+5z^3+14z^4+42z^5+132z^6+429z^7+1430z^8+4862z^9,$

\noindent $\bullet \quad z^2+3z^3+9z^4+28z^5+90z^6+297z^7+1001z^8+3432z^9,$

\noindent $\bullet \quad z^3+4z^4+14z^5+48z^6+165z^7+572z^8+2002z^9,$

\noindent which correspond to some shifts of  \oeis{A000108}, \oeis{A000245}, \oeis{A002057}, and \oeis{A000344} in \cite{Sloa}.


\section{Partial \L{}ukasiewicz paths constrained by height}

In this section, we count partial \L{}ukasiewicz paths of height at most $t$, for a  given $t\geq 0$. We introduce the notations $f_k^t$, $g_k^t$, $h_k^t$ for $0\leq k\leq t$,  $F^t(u)$, $G^t(u)$ and $H^t(u)$, which are the counterparts of $f_k$, $g_k$, $h_k$, $F(u)$, $G(u)$ and $H(u)$ for partial \L{}ukasiewicz paths of height at most $t\geq 0$. Considering the state diagram where each layer consists of only  $t+1$ states, we deduce a system of equations of the following form:

$$\left[\begin{array}{cccccccccc}
-1 & 0 & 0 & 0 & 0 & 0 & 0 & 0 & 0 & \cdots
\\
 0 & -1 & 0 & z  & z  & z  & 0 & 0 & 0 &\cdots
\\
 z  & z  & z -1 & 0 & 0 & 0 & 0 & 0 & 0 &\cdots
\\
 z  & z  & z  & -1 & 0 & 0 & 0 & 0 & 0 &\cdots
\\
 0 & 0 & 0 & 0 & -1 & 0 & z  & z  & z  &\cdots
\\
 0 & 0 & 0 & z  & z  & z -1 & 0 & 0 & 0 &\cdots
\\
 z  & z  & z  & z  & z  & z  & -1 & 0 & 0 &\cdots
\\
 0 & 0 & 0 & 0 & 0 & 0 & 0 & -1 & 0 &\cdots
\\
 0 & 0 & 0 & 0 & 0 & 0 & z  & z  & z -1 &\cdots
  \\
 \vdots & \vdots & \vdots & \vdots & \vdots & \vdots & \vdots  & \vdots  & \vdots &
\end{array}\right]\cdot \left[\begin{array}{c}
f_0^t
\\
g_0^t
\\
 h_0^t
\\
 f_1^t
\\
 g_1^t
\\
 h_1^t
\\
 f_2^t
\\
 g_2^t
\\
 h_2^t
 \\
 \vdots
\end{array}\right]
=\left[\begin{array}{c}
-1
\\
0
\\
 0
\\
 0
\\
 0
\\
0
\\
 0
\\
0
\\
0\\
\vdots\\
\end{array}\right]
$$

For a given height $t\geq 0$, the previous matrix (denoted $A_t$) is square with  $3(t+1)$ rows and  columns. Using classical properties of the determinant, we can prove that  $D_t=\det(A_t)$ satisfies
$$D_{t+2}+D_{t+1}+ zD_t=0 $$ anchored with $D_0 = z - 1$ and $D_1 = 1 - 2z$. Then we deduce, $$D_t=\frac{z(-1)^{t+2} \left(1-\sqrt{1-4 z}\right)^{t+2}}{2^{t+1}\sqrt{1-4 z}  \left(1+\sqrt{1-4 z}\right)}+\frac{z(-1)^{t+1} \left(1+\sqrt{1-4 z}\right)^{t+2}}{2^{t+1}\sqrt{1-4 z}  \left(1-\sqrt{1-4 z}\right)},$$ which corresponds to  $$D_t=(-1)^{t+1}\cdot F_t,$$ where $F_t$ is the Fibonacci polynomial (see \cite{Krew})  $$F_t=1-{t+1\choose 1}z+{t\choose 2}z^2-{t-1\choose 3}z^3+\ldots.$$

\noindent For instance, we have $D_3=F_3=1-4z+3z^2$ and $D_4=-F_4=-1+5z-6z^2+z^3$.

Resolving the system using the Cramer rule, we obtain for $0\leq k\leq t$

\begin{equation}
    f^{t}_k=\frac{N^t_{3k+1}}{D_t}, \quad g^{t}_k=\frac{N^t_{3k+2}}{D_t}, \quad h^{t}_k=\frac{N^t_{3k+3}}{D_t},
\label{eqsol}
\end{equation}
where $N^t_{k}$  is the determinant of the matrix $A_{t}(k)$ obtained from $A_{t}$ by replacing the $k$-th column with the vector ${}^T(-1,0,0,0,\ldots, 0,0)$.

As we have done for $D_t$, it is easy to prove that  $N^t_{k}$ satisfies the two recurrence relations:
$$N^{t+2}_k+N^{t+1}_{k}+zN^{t}_k=0, \quad 1\leq k,~ 1+\Big\lceil\frac{k}{3}\Big\rceil\leq t,  \mbox{ and }$$
$$N^{t+1}_{k+3}=-N^{t}_{k},\quad  4\leq k,~ 1\leq t.$$

Calculating the initial terms $N_{k}^{t}$, for $(t,k)\in\{0,1,2\}\times \{1,2,3\}$ and  for $(t,k)\in\{1,2,3\}\times\{4,5,6\}$, we can easily obtain  a closed form for all $N_{3k+i}^{t}$, $1\leq i\leq 3$, and then for all $f_k^{t}$, $g_k^{t}$, $h_k^{t}$, $0\leq k\leq t$. See Table~1 for several exact values of $N_k^t$.

\begin{table}
    \centering
   \begin{tabular}{cccccc}
k\slash t  & 0 & 1 & 2 & 3&4\\
\hline
1 & $z-1$ &$1-2z$ & $-(z^2-3z+1)$ &$3z^2-4z+1$&$z^3-6z^2+5z-1$ \\
2 & 0 &$z^2$ & $-z^2$&$z^2(1-z)$ &$-z^2(1-2z)$\\
3 & $-z$ &$z(1-z)$ &$-z(1-2z)$ &$z(z^2 -3z+1)$&$-z(3z^2-4z+1)$\\
4 &    &\cellcolor{lightlightgray}$z(1-z)$ &\cellcolor{llightgray}$-z(1-2z)$ &\cellcolor{lightlightgray} $z(z^2 -3z+1)$&\cellcolor{llightgray}$-z(3z^2-4z+1)$\\
5 &    &\cellcolor{lightlightgray}0 &\cellcolor{llightgray}$-z^2$ &\cellcolor{lightlightgray}$z^2$&\cellcolor{llightgray}$-z^2(1-z)$\\
6 &    & \cellcolor{lightlightgray} $z^2$ & \cellcolor{llightgray}$-z^2$& \cellcolor{lightlightgray}$z^2(1-z)$&\cellcolor{llightgray}$-z^2(1-2z)$\\
7 &    &      &\cellcolor{lightlightgray}$-z(1-z)$ &\cellcolor{llightgray}$z(1-2z)$&\cellcolor{lightlightgray}$-z(z^2-3z+1)$ \\
8 &    &      &\cellcolor{lightlightgray}0 & \cellcolor{llightgray}$z^2$&\cellcolor{lightlightgray}$-z^2$\\
9 &    &      &\cellcolor{lightlightgray}$-z^2$ & \cellcolor{llightgray}$z^2$&\cellcolor{lightlightgray}$-z^2(1-z)$\\
10 &    &      & &$\cellcolor{lightlightgray}z(1-z)$& \cellcolor{llightgray}$-z(1-2z)$\\
11 &    &      & &\cellcolor{lightlightgray} 0&\cellcolor{llightgray}$-z^2$\\
12 &    &      & &\cellcolor{lightlightgray} $z^2$&\cellcolor{llightgray}$-z^2$\\
13 &    &      & & &\cellcolor{lightlightgray}\ldots\\
\end{tabular}
\caption{The first values of $N_k^t$ for $0\leq t \leq 4$ and $1\leq k\leq 12$.}
\end{table}

In particular, we have $N_1^t=D_t$ (see above for a closed form),

\begin{equation*}
	N_2^t=\frac{z^{2} \left(-\frac{2 z}{1+\sqrt{1-4 z}}\right)^{t}}{\sqrt{1-4 z}}-\frac{z^{2} \left(-\frac{2 z}{-\sqrt{1-4 z}+1}\right)^{t}}{\sqrt{1-4 z}},
\end{equation*}
$$N_3^t=-\frac{z^{2} \left(-1+\sqrt{1-4 z}\right) \left(-\frac{2 z}{1+\sqrt{1-4 z}}\right)^{t}}{\sqrt{1-4 z}  \left(1+\sqrt{1-4 z}\right)}-\frac{z^{2} \left(1+\sqrt{1-4 z}\right) \left(-\frac{2 z}{-\sqrt{1-4 z}+1}\right)^{t}}{\sqrt{1-4 z}  \left(-\sqrt{1-4 z}+1\right)},
$$
and $N_4^t=N_3^t$, $N_5^t=-N_2^{t-1}$, $N_6^t=N_2^t$ for $t\geq 1$. Using (\ref{eqsol}), we can deduce easily the close  forms of $f_k^t$, $g_k^t$, $h_k^t$ for $0\leq k\leq 1$ and $k\leq t$. Using the recurrence relations for $N_k^t$, we deduce closed forms for the other cases.

\begin{thm} For $2\leq k\leq t$, we have

\begin{align}f_k&=[ u^k] F(u)=\frac{N_3^{t-k+1}}{D_t}(-1)^{k-1},\\
g_k&=[ u^k] G(u)=\frac{N_2^{t-k}}{D_t}(-1)^{k}, \mbox{ and }\\
h_k&=[ u^k] H(u)=\frac{N_2^{t-k+1}}{D_t}(-1)^{k-1}.\end{align}
\end{thm}
 The first values of the series expansion of $f_2$ for $t=2,3,4$ are:

 \noindent $\bullet \quad z+2z^2+5z^3+13z^4+34z^5+89z^6+233z^7+610z^8+1597z^9,$

 \noindent $\bullet \quad z+2z^2+5z^3+14z^4+41z^5+122z^6+365z^7+1094z^8+3281z^9,$

\noindent $\bullet \quad z+2z^2+5z^3+14z^4+42z^5+131z^6+417z^7+1341z^8+4334z^9,$
which correspond to the sequences \oeis{A001519},  \oeis{A007051}, \oeis{A080937} in \cite{Sloa}, that count Dyck paths of semilength $n$ of height at most $t+1$.
\begin{thm}

The  generating function $[ u^k] \mathit{Total}_t(z,u)$ for the  number of partial \L{}ukasiewicz paths  of height at most $t\geq 0$, ending at height $k\geq 1$, is given by

$$(-1)^{k-1}\frac{N_3^{t-k+1}-N_2^{t-k}+N_2^{t-k+1}}{D_t}.$$


Moreover, we have for $k=0$: $$[ u^0] \mathit{Total}_t(z,u)=\frac{D_t+N_2^t+N_3^t}{D_t}.$$

The generating function for the total number of partial \L{}ukasiewicz paths  of height at most  $t\geq 0$ is given by $$ \mathit{Total}_t(z,1)=(-1)^{t+1}\cdot D_t^{-1}=F_t^{-1}.$$
\end{thm}

The first values of the series expansion of $ \mathit{Total}_t(z,1)$ for $t=0,1,2,3,4$ are respectively

\noindent $\bullet \quad 1+z+z^2+z^3+z^4+z^5+z^6+z^7+z^8+z^9,$

\noindent $\bullet \quad 1+2z+4z^2+8z^3+16z^4+32z^5+64z^6+128z^7+256z^8+512z^9,$

\noindent $\bullet \quad 1+3z+8z^2+21z^3+55z^4+144z^5+377z^6+987z^7+2584z^8+6765z^9,$

\noindent $\bullet \quad 1+4z+13z^2+40z^3+121z^4+364z^5+1093z^6+3280z^7+9841z^8+29524z^9,$

\noindent $\bullet \quad 1+5z+19z^2+66z^3+221z^4+728z^5+2380z^6+7753z^7+25213z^8+81927z^9,$ which corresponds to   \oeis{A000012}, \oeis{A000079}, \oeis{A001906}, \oeis{A003462}, \oeis{A005021} in \cite{Sloa}.
\medskip

Using \cite{Krew}, $\left[ z^n\right] \mathit{Total}_t(z,1)$ counts also the paths of length $2n+1+t$ in $\Bbb{N}^2$ starting at the origin, ending at $(n+t+1,n)$, consisting of steps $(0,1)$, and $(1,0)$, and such that all its points $(x,y)$ satisfy $x-t-1\leq y\leq x$. It would be interesting to exhibit a constructive bijection between these paths and partial \L{}ukasiewicz paths of height at most  $t\geq 0$.


\section{Partial \L{}ukasiewicz paths constrained by height from right-to-left}

In this section, we count partial \L{}ukasiewicz paths from right-to-left of height at most $t$, for a  given $t\geq 0$. The notations $f_k^t$, $g_k^t$, $h_k^t$ for $0\leq k\leq t$,  $F^t(u)$, $G^t(u)$ and $H^t(u)$, are in the same way as for Section 4. We deduce a system of equations of the following form:

$$\left[\begin{array}{cccccccccc}
-1 & 0 & 0 & 0 & 0 & 0 & 0 & 0 & 0 & \cdots
\\
 0 & -1 & 0 & z  & z  & z  & z & z & z &\cdots
\\
 z  & z  & z -1 & 0 & 0 & 0 & 0 & 0 & 0 &\cdots
\\
 z  & z  & z  & -1 & 0 & 0 & 0 & 0 & 0 &\cdots
\\
 0 & 0 & 0 & 0 & -1 & 0 & z  & z  & z  &\cdots
\\
 0 & 0 & 0 & z  & z  & z -1 & 0 & 0 & 0 &\cdots
\\
 0  & 0  & 0  & z  & z  & z  & -1 & 0 & 0 &\cdots
\\
 0 & 0 & 0 & 0 & 0 & 0 & 0 & -1 & 0 &\cdots
\\
 0 & 0 & 0 & 0 & 0 & 0 & z  & z  & z -1 &\cdots
  \\
 \vdots & \vdots & \vdots & \vdots & \vdots & \vdots & \vdots  & \vdots  & \vdots &
\end{array}\right]\cdot \left[\begin{array}{c}
f_0^t
\\
g_0^t
\\
 h_0^t
\\
 f_1^t
\\
 g_1^t
\\
 h_1^t
\\
 f_2^t
\\
 g_2^t
\\
 h_2^t
 \\
 \vdots
\end{array}\right]
=\left[\begin{array}{c}
-1
\\
0
\\
 0
\\
 0
\\
 0
\\
0
\\
 0
\\
0
\\
0\\
\vdots\\
\end{array}\right]
$$

For a given height $t\geq 0$, the previous matrix (denoted $A'_t$) is square with  $3(t+1)$ rows and  columns. Using classical properties of the determinant, we can prove that  $\det(A'_t)$ equals $D_t=\det(A_t)$ for all $t\geq 0$.

Resolving the system using the Cramer rule, we obtain for $0\leq k\leq t$

\begin{equation}
    f^{t}_k=\frac{N^t_{3k+1}}{D_t}, \quad g^{t}_k=\frac{N^t_{3k+2}}{D_t}, \quad h^{t}_k=\frac{N^t_{3k+3}}{D_t},
\label{eqsol}
\end{equation}
where $N^t_{k}$  is the determinant of the matrix $A_{t}(k)$ obtained from $A_{t}$ by replacing the $k$-th column with the vector ${}^T(-1,0,0,0,\ldots, 0,0)$.

As we have done for $D_t$, it is easy to prove that  $N^t_{k}$ satisfies the two recurrence relations:
$$N^{t+2}_k+N^{t+1}_{k}+zN^{t}_k=0, \quad 1\leq k\leq 3,~2\leq t  \mbox{ and }$$
$$N^{t+1}_{k}=-zN^{t}_{k-3},\quad  4\leq k,~ 0\leq t.,$$
where for all $(t,k)\in \Bbb{N}\times \{1,2,3\}$, $N_k^t$ is the same as in Section 4 .

\begin{thm} For $0\leq k\leq t$, we have

\begin{align}f_k&=[ u^k] F(u)=\frac{N_1^{t-k}}{D_t}(-1)^{k}z^k,\\
g_k&=[ u^k] G(u)=\frac{N_2^{t-k}}{D_t}(-1)^{k}z^k, \mbox{ and }\\
h_k&=[ u^k] H(u)=\frac{N_3^{t-k}}{D_t}(-1)^{k}z^k.\end{align}
\end{thm}

\begin{thm}

The  generating function $[ u^k] \mathit{Total}_t(z,u)$ for the  number of partial \L{}ukasiewicz paths  (from right to left) of height at most $t\geq 0$, ending at height $k\geq 0$, is given by
$$(-1)^{k}z^k\frac{N_1^{t-k}+N_2^{t-k}+N_3^{t-k}}{D_t}.$$


The generating function for the total number of partial \L{}ukasiewicz paths (from right to left)  of height at most  $t\geq 2$ is given by
$$ \mathit{Total}_t(z,1)=\frac{D_{t-2}}{D_t},$$
with  $ \mathit{Total}_0(z,1)=\frac{1}{1-z},$ and  $ \mathit{Total}_1(z,1)=\frac{1}{1-2z}.$
\end{thm}

The first values of the series expansion of $ \mathit{Total}_t(z,1)$ for $t=0,1,2,3$ are respectively

\noindent $\bullet \quad 1+z+z^2+z^3+z^4+z^5+z^6+z^7+z^8+z^9,$

\noindent $\bullet \quad 1+2z+4z^2+8z^3+16z^4+32z^5+64z^6+128z^7+256z^8+512z^9,$

\noindent $\bullet \quad 1+2z+5z^2+13z^3+34z^4+89z^5+233z^6+610z^7+1597z^8+4181z^9,$

\noindent $\bullet \quad 1+2z+5z^2+14z^3+41z^4+122z^5+365z^6+1094z^7+3281z^8+9842z^9,$

 which corresponds to   \oeis{A000012}, \oeis{A000079}, \oeis{A001519}, \oeis{A007051} in \cite{Sloa}.
\medskip

\section{The average height of \L{}ukasiewicz paths}

\subsection{The left-to-right model}

We can simplify the expressions given in the previous section using the (magic) substitution $z=\frac{u}{(1+u)^2}$, first used in \cite{deBKR}. Then we find
\begin{gather*}
	D_t=(-1)^{t+3}\frac{1-u^{t+3}}{(1-u)(1+u)^{t+2}},\\
	N_2^t=\frac{(-1)^{t+1}u^2(1-u^t)}{(1+u)^{t+3}(1-u)},\\
	N_3^t= \frac{(-1)^{t+1}u(1-u^{t+2})}{(1-u)(1+u)^{t+3}}.	
\end{gather*}
We start with paths that return to the $x$-axis, as the formula is (slightly) different. We
have
\begin{equation*}
	\mathscr{L}^{[\le t]}=\frac{D_t+N_2^t+N_3^t}{D_t}=\frac{(1+u)-u^{t+2}(1+u^2)     }{1-u^{t+3}}
\end{equation*}
for the generating function of paths bounded by $t$. From this we see the limit for $t\to \infty$
$$\mathscr{L}^{[\le \infty]}=1+u=\frac{1-\sqrt{1-4z}}{2z},$$
as expected. From this
\begin{equation*}
	\mathscr{L}^{[>t]}=\mathscr{L}^{[\le \infty]}-\mathscr{L}^{[\le t]}=\frac{ u^{t+2}(1-u^2) }{1-u^{t+3}}.
\end{equation*}
We refer to \cite{HPW} where a similar instance is worked out with an extensive amount of detail. For the average height we have
(before normalizing by the Catalan numbers) to compute
\begin{equation*}	\sum_{t\ge0}\mathscr{L}^{[>t]}=\frac{1-u^2}{u}\sum_{t\ge3}\frac{ u^{t} }{1-u^{t}}.
\end{equation*}
The goal is to find the local behaviour of $u\sim1$ since it translates to the  local behaviour of $z\sim\frac14$. To find this, we set
$u=e^{-t}$ and use the Mellin transform; we don't need to do the actual computation, since we just cite the result from \cite{HPW}. First, the factor is simple;
\begin{equation*}
	\frac{1-u^2}{u}\sim 2(1-u).
\end{equation*}
Since we only compute the leading term of the asymptotic expansion, we use
\begin{equation*}
	\sum_{t\ge1}\frac{ u^{t} }{1-u^{t}}\sim -\frac{\log(1-u)}{1-u},
\end{equation*}
found e.g. in \cite{HPW}, so in total
\begin{equation*}
	-2\log(1-u)\sim-2\log(2\sqrt{1-4z})\sim-\log(1-4z).
\end{equation*}
Singularity analysis of generating functions \cite{FO} allows to translate this to the coefficients of $z^n$, with the result $\sim \frac{4^n}{n}$. For Catalan numbers, we have
the well-known
\begin{equation*}
	\frac1{n+1}\binom{2n}{n}\sim \frac{4^n}{n^{3/2}\sqrt{\pi}}.
\end{equation*}
The quotient is then the average height (leading term only):
\begin{equation*}
	\frac{\frac{4^n}{n}}{\frac{4^n}{n^{3/2}\sqrt{\pi}}}=\sqrt{\pi n}.
\end{equation*}

For path ending at $k\ge1$ we have to compute
\begin{equation*}
	(-1)^{k-1}\frac{N_3^{t-k+1}-N_2^{t-k}+N_2^{t-k+1}}{D_t}=u(1+u)^{k}\frac{1 - u^{t-k+1}}{1-u^{t+3}}.
\end{equation*}
The limit for $t\to \infty$ (no boundary) is $u(1+u)^{k}$ and the difference is
\begin{equation*}
	u(1+u)^{k}-u(1+u)^{k}\frac{1 - u^{t-k+1}}{1-u^{t+3}}=u^{-k-1}(1+u)^k(1-u^{k+2})\frac{u^{t+3}}{1-u^{t+3}}.
\end{equation*}
Locally we have ($u\sim1$)
\begin{equation*}
	u^{-k-1}(1+u)^k(1-u^{k+2})\sim 2^k(k+2)(1-u).
\end{equation*}
For the average height (leading term only, before normalization) we need to compute
\begin{equation*}
	2^k(k+2)(1-u)\sum_{t\ge1}\frac{u^{t}}{1-u^{t}}\sim-2^k(k+2)\log(1-u),
\end{equation*}
and
\begin{equation*}
	[z^n]-2^k(k+2)\log(1-u)\sim2^{k-1}(k+2)\frac{4^n}{n}.
\end{equation*}
For the total numbers, we have
\begin{equation*}
	[z^n]u(1+u)^k=\binom{2n-1+k}{n-1}-\binom{2n-1+k}{n-2}\sim \frac{4^n}{\sqrt{\pi}n^{3/2} }(k+2)2^{k-1},
\end{equation*}
and the average height ($k$ fixed, $n\to\infty$) is asymptotic to
\begin{equation*}
	\frac{2^{k-1}(k+2)\frac{4^n}{n}}{\frac{4^n}{\sqrt{\pi}n^{3/2} }(k+2)2^{k-1}}=\sqrt{\pi n},
\end{equation*}
as before.

To compute the average height of paths with unspecified endpoint makes no sense in this model since the number of such paths of length $n$ is infinite.

\subsection{The right-to left model}

We have to analyze
$$(-1)^{k}z^k\frac{N_1^{t-k}+N_2^{t-k}+N_3^{t-k}}{D_t}.$$
It translates into
\begin{equation*}
\mathscr{R}^{\le t}=\frac{u^k}{(1+u)^{k-1}}\frac{ 1	-u^{t+2-k}	}{1-u^{t+3}}.
\end{equation*}
Furthermore
\begin{equation*}
	\mathscr{R}^{\le \infty}=\frac{u^k}{(1+u)^{k-1}}
\end{equation*}
and
\begin{equation*}
\mathscr{R}^{>t}=	\mathscr{R}^{\le \infty}-\mathscr{R}^{\le t}=
\frac{u^{t+2}(1 -u^{k+1})}{(1+u)^{k-1}(1-u^{t+3})}.
\end{equation*}
For the average, we must compute
\begin{equation*}
	\frac{(k+1)(1 -u)}{2^{k-1}}	\sum_{t\ge1}\frac{u^{t}}{1-u^{t}}
\end{equation*}
where we took liberties about two missing terms, which does not influence the main term of the average height.

As before, we get the asymptotic behaviour
\begin{equation*}
	\frac{(k+1)}{2^{k}}	\frac{4^n}{n}.
\end{equation*}
For the total number we get
\begin{equation*}
\frac{k+1}{2^k}\frac{4^n}{\sqrt{\pi}n^{3/2}},
\end{equation*}
and the quotient is again $\sqrt{\pi n}$.

Now we move to the \L{}ukasiewicz paths with unspecified end and have to consider
\begin{equation*}
\frac{D_{t-2}}{D_t}
\end{equation*}
which is also
\begin{equation*}
\frac{(1+u)^3(1-u)}{u^2}\frac{u^{t+3}}{1-u^{t+3}}.
\end{equation*}
The fact that the formula is different for small values of $t$ and that we start the sum at 1 and not at 3, does not change the main term.
We get
\begin{equation*}
\frac{(1+u)^3(1-u)}{u^2}\sum_{t\ge1}\frac{u^{t}}{1-u^{t}}\sim -8\log(1-u),
\end{equation*}
and the coefficient of $z^n$ in it is asymptotic to
\begin{equation*}
\frac{4^{n+1}}{n};
\end{equation*}
for the total number of objects we get the asymptotic formula
\begin{equation*}
\frac{4^{n+1}}{\sqrt\pi  n^{3/2}}
\end{equation*}
and the quotient (average height) is again asymptotic to $\sqrt{\pi n}$.

\section{Partial alternate \L{}ukasiewicz paths}

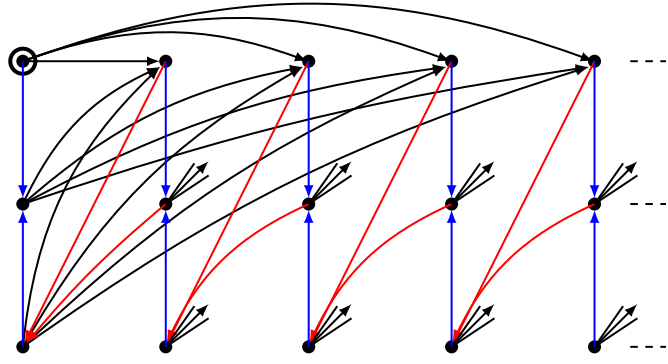
\begin{figure}[H]
\begin{center}\quad\scalebox{0.95}{\begin{tikzpicture}[ultra thick]
 \pointCercle{0cm}{4cm}
 \point{2cm}{4cm}\point{4cm}{4cm}\point{6cm}{4cm}\point{8cm}{4cm}
\point{0cm}{2cm}
\point{2cm}{2cm}\point{4cm}{2cm}\point{6cm}{2cm}\point{8cm}{2cm}
\point{0cm}{0cm}\point{2cm}{0cm}\point{4cm}{0cm}\point{6cm}{0cm}\point{8cm}{0cm}
  \draw[->,>=latex,black,thick] (0,4)--(1.92,4);
    \draw[black,dashed,thick] (8.5,4)--(9,4);
    \draw[black,dashed,thick] (8.5,2)--(9,2);
    \draw[black,dashed,thick] (8.5,0)--(9,0);
  \draw[->,>=latex,thick] (0,4) to[bend left=20] (3.99,4);
  \draw[->,>=latex,thick] (0,4) to[bend left=20] (5.99,4);
  \draw[->,>=latex,thick] (0,4) to[bend left=20] (7.99,4);
%
  \draw[->,>=latex,thick] (0,2) to[bend left=20] (1.92,3.92);
    \draw[->,>=latex,thick] (0,2) to[bend left=15] (3.92,3.92);
     \draw[->,>=latex,thick] (0,2) to[bend left=10] (5.92,3.92);
     \draw[->,>=latex,thick] (0,2) to[bend left=5] (7.92,3.92);
     \draw[->,>=latex,thick] (0,0) to[bend left=20] (1.92,3.92);
    \draw[->,>=latex,thick] (0,0) to[bend left=15] (3.92,3.92);
     \draw[->,>=latex,thick] (0,0) to[bend left=10] (5.92,3.92);
     \draw[->,>=latex,thick] (0,0) to[bend left=9] (7.92,3.92);
     \draw[->,>=latex,red,thick] (2,2) to[bend right=5] (0.05,0.05);
      \draw[->,>=latex,red,thick] (4,2) to[bend right=20] (2.01,0.05);
     \draw[->,>=latex,red,thick] (6,2) to[bend right=20] (4.01,0.05);
     \draw[->,>=latex,red,thick] (8,2) to[bend right=20] (6.01,0.05);
     \draw[->,>=latex,red,thick] (2,4)-- (0.02,0.05);
     \draw[->,>=latex,red,thick] (4,4)-- (2.02,0.05);
     \draw[->,>=latex,red,thick] (6,4)-- (4.02,0.05);
     \draw[->,>=latex,red,thick] (8,4)-- (6.02,0.05);
      \draw[thick] (2,0)-- (2.4,0.57);
      \draw[->,>=latex,thick] (2,0)-- (2.6,0.6);
       \draw[thick] (2,0)-- (2.6,0.4);
       \draw[thick] (2,2)-- (2.4,2.57);
      \draw[->,>=latex,thick] (2,2)-- (2.6,2.6);
       \draw[thick] (2,2)-- (2.6,2.4);
        \draw[thick] (4,0)-- (4.4,0.57);
      \draw[->,>=latex,thick] (4,0)-- (4.6,0.6);
       \draw[thick] (4,0)-- (4.6,0.4);
       \draw[thick] (4,2)-- (4.4,2.57);
      \draw[->,>=latex,thick] (4,2)-- (4.6,2.6);
       \draw[thick] (4,2)-- (4.6,2.4);
        \draw[thick] (6,0)-- (6.4,0.57);
      \draw[->,>=latex,thick] (6,0)-- (6.6,0.6);
       \draw[thick] (6,0)-- (6.6,0.4);
       \draw[thick] (6,2)-- (6.4,2.57);
      \draw[->,>=latex,thick] (6,2)-- (6.6,2.6);
       \draw[thick] (6,2)-- (6.6,2.4);
        \draw[thick] (8,0)-- (8.4,0.57);
      \draw[->,>=latex,thick] (8,0)-- (8.6,0.6);
       \draw[thick] (8,0)-- (8.6,0.4);
       %
       %
       %
      \draw[thick] (8,2)-- (8.4,2.57);
     \draw[->,>=latex,thick] (8,2)-- (8.6,2.6);
     \draw[thick] (8,2)-- (8.6,2.4);
       %
\draw[->,>=latex,blue, thick] (0,0)--(0,1.92);
\draw[->,>=latex,blue, thick] (2,0)--(2,1.92);
\draw[->,>=latex,blue, thick] (4,0)--(4,1.92);
\draw[->,>=latex,blue, thick] (6,0)--(6,1.92);
\draw[->,>=latex,blue, thick] (8,0)--(8,1.92);
\draw[->,>=latex,blue, thick] (0,4)--(0,2.08);
\draw[->,>=latex,blue, thick] (2,4)--(2,2.08);
\draw[->,>=latex,blue, thick] (4,4)--(4,2.08);
\draw[->,>=latex,blue, thick] (6,4)--(6,2.08);
\draw[->,>=latex,blue, thick] (8,4)--(8,2.08);
 %
%
 \end{tikzpicture}}
\end{center}
\caption{The state diagram for partial alternate \L{}ukasiewicz paths. Black (resp.\ red, resp.\ blue) arrows corresponds to up-steps (resp.\ down-steps, resp.\ horizontal steps).}
\label{fig4}\end{figure}

In this section, we study  prefixes of alternate \L{}ukasiewicz  paths, i.e., \L{}ukasiewicz  paths that do not contain two consecutive steps with the same direction (or equivalently,  walks in the state diagram of Figure~\ref{fig3} without two consecutive arrows of the same color). We refer to Figure~\ref{fig4} for the state diagram associated to these paths. The notations $f_k$, $g_k$ and $h_k$, $k\geq 0$, are already defined previously, but have a different meaning in this section.

\begin{equation}
\begin{array}{ll}
f_0&=1,\quad \mbox{and}\quad f_k=zf_0+z\sum\limits_{\ell=0}^{k-1} g_\ell +z\sum\limits_{\ell=0}^{k-1} h_\ell, \quad k\geq 1,\\
g_k&=z f_{k+1}+zh_{k+1},\quad k\geq 0,\\
h_k&=zf_k+zg_k,\quad k\geq 0.
\end{array}
\label{equalt}
\end{equation}

Considering generating functions
$$F(u,z)=\sum\limits_{k\geq 0} u^kf_k(z), \quad G(u,z)=\sum\limits_{k\geq 0} u^kg_k(z), \mbox{ and } H(u,z)=\sum\limits_{k\geq 0} u^kh_k(z),$$
and summing the recursions in (\ref{equalt}), we obtain:
\begin{align*}
F(u)&=1+z\sum\limits_{k\geq 1}u^k\biggl(1+\sum\limits_{\ell=0}^{k-1} g_\ell +\sum\limits_{\ell=0}^{k-1} h_\ell\biggr)\\
&=1+\frac{zu}{1-u}+z\sum\limits_{k\geq 0}\frac{u^{k+1}}{1-u}g_k+z\sum\limits_{k\geq 0}\frac{u^{k+1}}{1-u}h_k\\
&=1+\frac{zu}{1-u}(1+G(u)+H(u)),\\
G(u)&=\frac{z}{u}(F(u)+H(u)-1-H(0)),\\
H(u)&=z(F(u)+G(u)).
\end{align*}

Solving these functional equations, we deduce
\begin{gather*}
F(u)=\frac{uz^2(1+z)H(0)+2 u  z^{3}-u^{2} z +u^{2}+z^{2}-u}{u^{2} z^{2}+2 u  z^{3}+u^{2}+z^{2}-u},\\
 G(u)=-\frac{z \left(H(0)(u z^{2}+u-1)+2 u  z^{2} +z \right)}{u^{2} z^{2}+2 u  z^{3}+u^{2}+z^{2}-u},\\
H(u)=\frac{z \left(zH(0)(u  z-u+1) -u^{2} z +u^{2}-u \right)}{u^{2} z^{2}+2 u  z^{3}+u^{2}+z^{2}-u}.
\end{gather*}

Now we apply the kernel method on $H(u)$. We have
\begin{equation}
    H(u)=\frac{z \left(zH(0)(u  z-u+1) -u^{2} z +u^{2}-u \right)}{(1+z^2)(u-s_1)(u-s_2)}
\label{equH}
\end{equation}
with
$$s_1=\frac{1-2 z^{3}+\sqrt{4 z^{6}-4 z^{4}-4 z^{3}-4 z^{2}+1}}{2 z^{2}+2},$$

$$s_2=\frac{1-2 z^{3}-\sqrt{4 z^{6}-4 z^{4}-4 z^{3}-4 z^{2}+1}}{2 z^{2}+2}.$$

 In order to compute $H(0)$, it suffices to plug $u=s_2$ in the numerator of (\ref{equH}). Then, $H(0)$ satisfies $zH(0)(s_2  z-s_2+1) -s_2^{2} z +s_2^{2}-s_2=0$, which implies that
 $$H(0)=\frac{s_2}{z}.$$
 After this, and using $s_1s_2(1+z^2)=z^2$, we simplify  of both, numerators and denominators, in $F(u)$, $G(u)$, $H(u)$  by factorizing them with $(u-s_2)$.

 $$ F(u)=\frac{(1-z)u-(1+z^2)s_1}{(1+z^2)(u-s_1)}=-\frac{z -1}{z^{2}+1}-\frac{s_1 z \left(z +1\right)}{\left(u -s_1 \right) \left(z^{2}+1\right)},$$
 $$G(u)=-\frac{z^2s_1^{-1}+2z^3}{(1+z^2)(u-s_1)}=-\frac{z^2(2zs_1+1)}{(1+z^2)s_1(u-s_1)},$$
 $$ H(u)=\frac{z((1-z)u-1)}{(1+z^2)(u-s_1)}=\frac{\left(1-z\right) z}{z^{2}+1}-\frac{\left( zs_1 -s_1 +1\right) z}{\left(u -s_1 \right) \left(z^{2}+1\right)}.$$

Finally, we easily obtain:
\begin{align}
f_k&=[ u^k] F(u)=\frac{z(z+1)}{(1+z^2)s_1^k},\\
g_k&=[ u^k] G(u)=\frac{z^2(2zs_1+1)}{(1+z^2)s_1^{k+2}}, \quad \mbox{and}\\
h_k&=[ u^k] H(u)=\frac{z((z-1)s_1+1)}{(1+z^2)s_1^{k+1}}.
\end{align}

 Since the expansion of $s_1$ does not have pretty coefficients, we cannot expect this from our final answers.

\begin{thm} The bivariate generating function for the total number of partial alternate \L{}ukasiewicz paths of length $n$  with respect to the length $n$ and the height of the last point is given by
$$\mathit{Total}(z,u)=\frac {s_1^{2}{z}^{2}+s_1 u{z}^{2}+2 s_1 {z}^{3}+
s_1^{2}-s_1 u+zs_1+{z}^{2}}{ \left( {z}^{2}+1 \right)  ( -u+s_1 ) s_1}
$$ and we have
$$[ u^0] \mathit{Total}(z,u)={\frac {s_1^{2}{z}^{2}+2 s_1 {z}^{3}+s_1^{2}+ s_1z+{z}^{2}}{ ( {z}^{2}+1 ) s_1^{2}}},
$$
and for $k\geq 1$,
$$[ u^k] \mathit{Total}(z,u)={\frac {z ( 2s_1^{2}z+2s_1{z}^{2}+s_1+z) }{ ( {z}^{2}+1 ) s_1^{k+2}}}
.$$

\end{thm}
Here are examples of  the series expansions of $[ u^k]\mathit{Total}(z,u)$ for $k=0,1,2,3$ (leading terms):

\noindent $\bullet \quad 1+z+z^2+3z^3+5z^4+9z^5+19z^6+39z^7+81z^8+173z^9,$

\noindent $\bullet \quad z+3z^2+5z^3+11z^4+25z^5+53z^6+115z^7+255z^8+565z^9;$

\noindent $\bullet \quad z+3z^2+7z^3+19z^4+45z^5+105z^6+247z^7+575z^8+1333z^9,$

\noindent $\bullet \quad z+3z^2+ 9z^3+27z^4+69z^5+177z^6+443z^7+1087z^8+2645z^9,$

\noindent which do not appear in \cite{Sloa}. As a consequence, the first terms of the series expansion of the generating function for the number of alternate \L{}ukasiewicz paths are  $$1+z+z^2+3z^3+5z^4+9z^5+19z^6+39z^7+81z^8+173z^9.$$

Singularity analysis of generating function $[ u^0] \mathit{Total}(z,u)$ gives
$$[ z^n][ u^0]\mathit{Total}(z,u)\sim \frac{\sqrt{-6 a^{6}+4 a^{4}+3 a^{3}+2 a^{2}}  (a +1) 2^{n} \left(-a^{2}+a +1\right)^{n}}{\sqrt{\pi}  a^{2} \left(a^{2}+1\right) n^{\frac{3}{2}}},$$
with $$a=\frac{1}{3}-\frac{2 \cos  \Bigl(\frac{\arctan \left(\frac{15 \sqrt{111}}{487}\right)}{6}+\frac{\pi}{6}\Bigr)}{3}+\frac{2 \sin  \Bigl(\frac{\arctan \left(\frac{15 \sqrt{111}}{487}\right)}{6}+\frac{\pi}{6}\Bigr) \sqrt{3}}{3}.$$

The reason that this constant appears is this: For singularity analysis, one needs the solution closest to the origin of $4 z^{6}-4 z^{4}-4 z^{3}-4 z^{2}+1=0$,
which is $0.403031716762\dots$. Maple gives the curious explicit version if one asks for a `simplification.'

\end{document}